\documentstyle[12pt]{article}
\renewcommand{\baselinestretch}{1.3}

\makeatletter
\newcommand{\single}{\let\CS=\@currsize\renewcommand{\baselinestretch}{1}\tiny\CS}
\newcommand{\singles}{\let\CS=\@currsize\renewcommand{\baselinestretch}{1.3}\tiny\CS}
\newcommand{\oneanda}{\let\CS=\@currsize\renewcommand{\baselinestretch}{1.2}\tiny\CS}
\newcommand{\doubles}{\let\CS=\@currsize\renewcommand{\baselinestretch}{1.5}\tiny\CS}
\newcommand{\tree}{\let\CS=\@currsize\renewcommand{\baselinestretch}{1.5}\tiny\CS}
\newcommand{\four}{\let\CS=\@currsize\renewcommand{\baselinestretch}{2}\tiny\CS}

\newcommand{\ncom}{\newcommand}
\ncom{\bq}{\begin{equation}}
\ncom{\eq}{\end{equation}}
\ncom{\beqn}{\begin{eqnarray*}}
\ncom{\eeqn}{\end{eqnarray*}}
\ncom{\beq}{\begin{eqnarray}}

\ncom{\eeq}{\end{eqnarray}}
\ncom{\been}{\begin{enumerate}}
\ncom{\eeen}{\end{enumerate}}
\ncom{\nno}{\nonumber}
\ncom{\hs}{\mbox{\hspace{.25cm}}}
\ncom{\rar}{\rightarrow}
\ncom{\lrar}{\longrightarrow}
\ncom{\Rar}{\Rightarrow}
\ncom{\noin}{\noindent}
\newtheorem{thm}{Theorem}[section]
\newtheorem{lemma}[thm]{Lemma}
\newtheorem{cor}[thm]{Corollary}
\newtheorem{pro}[thm]{Proposition}
\newtheorem{example}[thm]{Example}
\newtheorem{remark}[thm]{Remark}

\ncom{\bt}{\begin{thm}}
\ncom{\et}{\end{thm}}
\ncom{\bl}{\begin{lemma}}
\ncom{\el}{\end{lemma}}
\ncom{\bco}{\begin{cor}}
\ncom{\eco}{\end{cor}}
\ncom{\bp}{\begin{pro}}
\ncom{\ep}{\end{pro}}
\ncom{\bex}{\begin{example}}
\ncom{\eex}{\end{example}}
\ncom{\brm}{\begin{remark}}
\ncom{\erm}{\end{remark}}
\ncom{\comx}{I\!\!\!\!C}
\ncom{\zee}{$Z\!\!\!\!Z$}
\ncom{\ze}{Z\!\!\!\!Z}
\ncom{\Q}{$I\!\!\!\!Q$}
\ncom{\p}{I\!\!P}
\ncom{\s}{(\!\!\times}
\ncom{\al}{\alpha}
\ncom{\be}{\beta}
\ncom{\f}{\frac}
\ncom{\ga}{\gamma}
\ncom{\bib}{\bibitem}
\ncom{\pf}{{\bf Proof: }}

\ncom{\sta}{\stackrel}

\ncom{\cA}{{\cal A}}
\ncom{\cG}{{\cal G}}
\ncom{\cI}{{\cal I}}
\ncom{\cO}{{\cal O}}
\ncom{\cV}{{\cal V}}
\ncom{\cW}{{\cal W}}
\ncom{\cK}{{\cal K}}
\ncom{\cM}{{\cal M}}
\ncom{\cL}{{\cal L}}
\ncom{\cZ}{{\cal Z}}
\title{Line bundles of type $(1,...,1,2,...,2,4,...,4)$ on Abelian Varieties}
\author{Jaya N.Iyer\\Institut de mathematiques, Case 247\\Univ. Paris -6, 4,Place Jussieu,
\\75252 Paris Cedex 05, France\\(email: iyer@math.jussieu.fr)}

\begin{document} 
\maketitle

\begin{abstract}We show  birationality of the morphism associated
  to line bundles $L$ of type $(1,...,1,2,...,2,4,...,4)$ on a generic
  $g-$dimensional abelian variety into its complete linear system such
  that $h^0(L)=2^g$. When $g=3$, we describe the image of the abelian
  threefold and from the
geometry of the moduli space $SU_C(2)$ in the linear system
  $|2\theta_C|$, we obtain analogous results in $\p H^0(L)$.
\end{abstract}

Mathematics Classification Number: 14C20, 14J17, 14J30, 14K10, 14K25.

\section{Introduction}
Let $L$ be an ample line bundle of type $\delta= (\delta_1,\delta_2,...,\delta_g)$ on a $g$-dimensional abelian variety $A$. Consider the associated rational
map $\phi_L:A\lrar \p H^0(A,L)$.

When $g=2$, Birkenhake, Lange and van Straten ( see [3]) have studied
line bundles of type $(1,4)$ on abelian surfaces. 
Suppose $L$ is an ample line bundle of type $(1,4)$ on an abelian
surface $A$. Then there is a cyclic covering $\pi:A\lrar B$ of degree
$4$ and a line bundle $M$ on $B$ such that $\pi^*M=L$. Let $X$ denote
the unique divisor in $|M|$ and put $Y=\pi^{-1}(X)$.
Their main theorem 
is
\bt
1) $\phi_L: A\lrar A'\subset \p ^3$ is birational onto a
singular octic $A'$ in $\p ^3$ if and only if $X$ and $Y$ do not admit
elliptic involutions compatible with the action of the Galois group of
$\pi$.

2)In the exceptional case $\phi_L:A\lrar A'\subset \p^3$ is a double
  covering of a singular quartic $A'$, which is birational to an
  elliptic scroll.

\et

Here we generalise this situation to higher dimensions and show

\bt
Suppose $L$ is an ample line bundle of type $\delta=(1,...,1,2,...,2,4,...,4)$
on a $g$-dimensional abelian variety $A$, $g\geq3$, such that $1$ and $4$ occur
equally often and atleast once in $\delta$. Then, for a generic pair $(A,L)$, the following holds. 

a) The associated morphism $\phi_L: A\lrar \p H^0(A,L)$ is birational
 onto its image.

b) When $g=3$, the image $\phi_L(A)$, can be described as follows,

there are $4$ curves $C_i$ on the image $\phi_L(A)$ such that the restricted morphism
 $\phi_L: \phi_L^{-1}(C_i)\lrar C_i\subset \phi_L(A)$ is of degree $2$.

\et

 Birkenhake et.al (see [3], Proposition 1.7, p.631) have shown the existence of the following commutative diagram

$
\begin{array}{ccccc}
A& \sta{\phi_L}{\lrar} & \phi_L(A) & \subset & \p^3=\p H^0(L)\\
\downarrow{\!\pi} & &\downarrow & &\downarrow{\!p}  \\
B& \sta{\phi_{M^2}} {\lrar} & \cK(B) & \subset & \p^3=\p H^0(M^2)\\
\end{array}
$

where $p(z_0:z_1:z_1:z_3)=(z_0^2:z_1^2:z_2^2:z_3^2)$ and the pair
$(B,M)$ is a principally polarized abelian surface. This diagram explains the geometry of
the image $\phi_L(A)$ from the geometry of the Kummer surface $\cK(B)$ and it 
 also gives the explicit equation of the surface $\phi_L(A)$ in $\p^3$.

Similarly, when $g\geq 3$ and the pair $(A,L)$ as in 1.2, we show that there is
 a commutative diagram:

$
\begin{array}{ccccc}
A& \sta{\phi_L}{\lrar} & \phi_L(A) & \subset & \p^{2^g-1}=\p H^0(L\\
\downarrow{\!\pi} & &\downarrow & &\downarrow{\!p}  \\
B & \sta{\phi_{M^2}} {\lrar} & \cK(B) & \subset & \p^{2^g-1}=\p H^0(M^2)\\
\end{array}
$

where $p(z_0:...:z_{2^g-1})=(z^2_0:...:z^2_{2^g-1})$ and $\pi$ is an
isogeny of degree $2^g$ and the pair $(B,M)$ is a principally polarized abelian variety.
This will explain the birationality of the map $\phi_L$ and the geometry of
the image $\phi_L(A)$, when $g=3$, as asserted in 1.2.
Since $deg(\phi_{M^2}\circ\pi)=2^{g+1}$ and from the birationality of $\phi_L$, it follows
that $deg(p|_{\phi_L(A)})= 2^{g+1}$. But since $degp=2^{2^g-1}$  the inverse image of the Kummer variety in $\p H^0(L)$ has components other than the
image $\phi_L(A)$. Hence the image $\phi_L(A)$ will be defined by forms other than
those coming from those forms which define the variety $\cK(B)$.
 
We study the situation when $g=3$, in detail.
Consider a pair $(A,L)$, with $L$ being an ample line bundle of 
type $(1,2,4)$ on an abelian threefold $A$. Consider an isogeny $A\lrar 
B= A/G$, where $G$ is a maximal isotropic subgroup of $K(L)$ of the
type $\frac{\ze}{2\ze}\times \frac{ \ze}{2\ze}\times \frac{\ze}{2\ze}$. Then $B$ is a principally polarized abelian threefold.
 If $B$ is isomorphic to the Jacobian variety of $C$, $J(C)$, where $C$ is a
smooth non-hyperelliptic curve of genus $3$, then the situation becomes
interesting because of the following results due to Narasimhan and Ramanan.

\bt
(See [12], Main Theorem, p.416) If $C$ is a non-hyperelliptic curve of genus $3$, then the moduli
space $SU_C(2)$ is isomorphic to a quartic hypersurface in $\p ^7$.
\et
( Here $\p^7= |2\theta|$, where $\theta$ is the canonical principal
polarization on the Jacobian $J(C)$ and $SU_C(2)$ is the moduli space
of rank $2$ semi-stable vector bundles with trivial determinant on the
curve $C$).

\bt
( See [11]) The Kummer variety $\cK$ is precisely the singular locus of $SU_C(2)$,
if $g(C)\geq 3$.
\et

The quartic hypersurface, $F=0$, is classically called the $\it{Coble\, quartic}$
 and is $\cG(2\theta)$-invariant in the linear system $|2\theta|$.
We identify the group of projective transformations, $H$, of order $8$, which
 acts on $\pi^{-1}\cK (C)$, (see 3.7). The $\cG(L)$-invariant octic hypersurface
$R$, given as $F(z^2_0:...:z^2_7)=0$ in $\p H^0(L)$, then contains the 
components $h(\phi_L(A)),h\in H$ in its singular locus.

Now we use the geometry of the moduli space $SU_C(2)$ in the linear system
$|2\theta|$, which has been extensively studied ( see [5], for instance),
 to get analogous results in $\p H^0(L)$.

We show
\bt
Consider a pair $(A,L)$, as above. Let $a\in K(L)$ be an element of
order $2$ such that $e^L(a,g)=-1$, for all $g\in G$, (here $e^L$ is
the Weil form on the group $K(L)$). Let $\p W_a$ be
an eigenspace in $\p H^0(L)$, for the action of $a$. Then there is a 
polarized abelian surface $(Z,N)$, $N$ is ample of type $(1,4)$ and a
commutative diagram

$
\begin{array}{cccccc}
Z& \sta{\phi_N}{\lrar} & \phi_N(Z) & \subset & \p H^0(N)\simeq & \p W_a \\
\downarrow{\!f} & &\downarrow & &\downarrow{\!q} & \downarrow{\!p} \\
P_a& \sta{\phi_{2\theta_a}} {\lrar} & \cK(P_a) & \subset &
\p H^0(2\theta_a)\simeq & \p V_a \\
\end{array}
$

Here $(P_a,\theta_a)$ is the Prym variety associated to the
$2$-sheeted unramified cover of the curve $C$, given by $\pi(a)$ and $\p V_a$ is the
eigenspace in $\p H^0(2\theta)$, for the action of $\pi(a)$.
The isomorphisms above are Heisenberg equivariant and the morphism
$q$ is given as $(r_0:r_1:r_2:r_3)\mapsto (r^2_0:r^2_1:r^2_2:r^2_3)$.
\et

We thus obtain the situation described by Birkenhake et.al in the case
$g=2$, nested in the case $g=3$. 

Moreover, the $\cG (N)$-invariant octic surface $\phi_N(Z)$ is mapped isomorphically onto
the $a^\perp/a(\simeq Heis(4))$-octic $R\cap \p W_a$ and we identify the set $\cap_{h\in H}h(
\phi_L(A))$ with the set of all pinch points and the coordinate points in $\phi_N(Z)$, occurring in each of
the eigenspace $\p W_a$, (see 5.6).
Finally, we make some remarks on the moduli space $\cA^{(1,2,4)}$.

$\it{Acknowledgements}$: We thank W.M.Oxbury and B.van Geemen for
making useful comments in an earlier version.
We are grateful to Christian Pauly for suggestions during revision.
We also thank the French Ministry of National Education, Research and
Technology, for their support.

$\bf{Notation:}$
Suppose $L$ is a symmetric line bundle i.e. $L\simeq i^*L$ for the involution
$i:A\lrar A$, $a\mapsto -a$.

The $\it{fixed\, group}$ of $L$ is $K(L)=\{a\in A: L\simeq
t_a^*L\},\,t_a:A\lrar A,x\mapsto a+x$.

The $\it{theta\, group}$ of $L$ is
$\cG(L)=\{(a,\phi):L\sta{\phi}{\simeq}t_a^*L\}$.

$K_1(\delta)=\frac{\ze}{d_1\ze}\times ...\times\frac{\ze}{d_g\ze}$, 
and $\widehat{K_1(\delta)} =Hom(K_1(\delta),\comx^*)$.

The $\it{Heisenberg\, group}$ of type $\delta$,
$Heis(\delta)=\comx^*\times K_1(\delta)\times \widehat{K_1(\delta)}$
and 
$V(\delta)=\{f: f:K_1(\delta)\lrar \comx\}$.

The $\it{Weil\, form}\, e^L:K(L)\times K(L)\lrar \comx^*$, is the
commutator map $(x,y)\mapsto x'y'x'^{-1}y'^{-1}$, for any lifts
$x',y'\in \cG(L)$ of $x,y\in K(L)$. 

For any $a\in K(L)$, $a^\perp= \{x\in K(L): e^L(a,x)=1\}$.

Consider the semi-direct product, $\cG(L)\s (i)$, of the theta group
associated to $L$ and the group generated by the involution $i$.
Let $\gamma\in \cG(L)\s (i)$ be an element of order 2.

$H^0(L)^{\pm}_\gamma = (\pm 1)-$eigenspace of $H^0(L)$ for the action of $\gamma$.

$h^0(L)^{\pm}_\gamma = dim H^0(L)^{\pm}_\gamma.$

$Q(V)=$ function field of a variety $V$.

\section{Birationality of the map $\phi_L$.}
 
Let $L$ be an ample line bundle of type $\delta=(1,..2,..,4)$ on a
$g$-dimensional abelian variety $A$. Here number of 2's= number of 4's in $\delta$.
Let $K(L)=\{a\in A: t^*_aL\simeq L\}$, where $t_a$ denotes translation by
$a$ on $A$. Choose a maximal isotropic subgroup $G$ of $K(L)$ w.r.t. the
Weil form $e^L$, containing $2K(L)$ and having only elements of order 2.
Then $G\simeq \frac{\ze}{2\ze}\times...\times \frac{\ze}{2\ze}$, $g$-times.
Consider the exact sequence
$$ 1\lrar \comx^*\lrar \cG(L)\lrar K(L)\lrar 0.$$
Let $G'$ be a lift of $G$ in $\cG(L)$.
Consider the isogeny $A\sta{\pi}{\lrar}B=A/G$. Then $L$ descends to a 
principal polarization $M$ on $B$. 
By Projection formula and using the fact that $\pi_*{\cal O} _A= \oplus_{\chi\in
\hat{G} }L_\chi$, where $L_\chi$ denotes the line bundle corresponding to
the character $\chi$, we deduce that 

$$ H^0(L)=\oplus_{\chi\in \hat{G}}H^0(M\otimes L_\chi).$$
Hence $\{s_\chi \in H^0(M\otimes L_\chi) : \chi\in \hat{G}\}$ is a basis
for the vector space $H^0(L)$ and since $M^2\otimes L^2_{\chi}\simeq M^2$,
$s^2_\chi=s_{\chi}\otimes s_\chi\in H^0(M^2) \forall \chi \in \hat{G}$.

Consider the homomorphism $\epsilon_2:\cG(L)\lrar \cG(L^2),
(x,\phi)\mapsto (x,\phi^{\otimes 2})$ and the inclusion $K(L)\subset K(L^2)$.

Then the subgroup $G\subset K(L^2)$ is isotropic for the
Weil form $e^{L^2}$.
Moreover, if $x\in K(L)$ and $g\in G$, then
$$ e^{L^2}(x,g)=e^L(x,g).e^L(x,g)=1.$$
Hence $\epsilon_2(\cG(L))\subset \cZ(\epsilon_2(G^{'}))$ and $\pi(K(L))\subset
K(M^2)$. ( Here $\cZ(\epsilon_2(G^{`}))= \{a\in
\cG(L^2): a.g^{'}=g^{'}.a, \forall g^{'}\in \epsilon _2(G^{'})\}$).

Now $\cG(M^2)= \cZ(\epsilon_2(G^{'}))/\epsilon _2(G^{'})$ and $H^0(M^2)=
H^0(L^2)^{G^{'}}$, where $H^0(L^2)^{G^{'}}$ denotes the vector subspace of
$\epsilon _2(G^{'})$-fixed sections of $H^0(L^2)$. For $g'\in G'$ and
$\chi\in \hat{G}$, $g'(s^2_\chi)= \chi^2(g).s^2_\chi=s^2_\chi$.
Hence $s^2_\chi\in H^0(L^2)^{G'}$, for all $\chi\in \hat{G}$.

We now show that
 $\{s^2_\chi : \chi\in \widehat{G}\} $ is a basis for $H^0(M^2)$, for a generic pair $(A,L)$..

In fact, we show that the homomorphism
$$\sum_{\chi\in \hat{G}}H^0(M\otimes L_\chi).H^0(M\otimes L_\chi)
\sta{\rho}{\lrar} H^0(M^2) ... (*)$$
is an isomorphism, for a generic pair $(A,L)$.

Consider the pair $(A,L)=(E_1\times ...\times E_r \times,
A_1\times...A_s, p_1^* L_1 \otimes...\otimes p_{r+s}^* L_{r+s})$, where $r$
is the number of 2's occurring in $\delta$,
 $E_1,...,E_r$ are elliptic curves with line bundles $L_i$ on $E_i$
of degree 2 and $A_j$ are simple abelian surfaces with line bundles $L_j$ on
$A_j$ of type $(1,4)$ ( by 1.1, $\phi_{L_j}(A_j)\subset |L_j|$ is an
octic surface).

In this case, one can easily see that the homomorphism 
$$S=Sym^2H^0(L_1)\otimes...\otimes Sym^2H^0(L_{r+s})\lrar
H^0(L_1^2)\otimes...\otimes H^0(L_{r+s}^2)= H^0(L_1^2\otimes
...\otimes L_{r+s}^2)$$
is injective.
Here, $(B,M)=(F_1,M_1)\times...\times (F_r,M_r)\times
(B_1,M'_1)\times...(B_s,M'_s)$, where $(F_j,M_j)$ are polarised
elliptic curves of degree 1 and $(B_j,M_j)$ are principally polarised
abelian surfaces. Also, the group $G$ is generated by elements of the
type $(e_1,..,e_r,a'_{r+1},..,a'_g)$, where each of $e_j$ and $a'_j$
are non-trivial 2 torsion elements of $E_j$ and $A_j$, respectively. 
Now it is easy to see that
$\sum_{\chi\in \hat{G}}H^0(M\otimes L_{\chi}).H^0(M\otimes
L_\chi)\subset S$ and $H^0(M^2)\subset H^0(L^2)$ and (*) is an isomorphism.

Hence, for a generic pair $(A,L)$ as above, (*) is an isomorphism.

As a consequence, we obtain the following
\bp
Consider a generic principally polarized abelian variety $(B',M')$ of 
dimension $g$. Let $H$ be a subgroup of 2-torsion points of $B'$, of
order g. Then the image of $H$ in $\cK(B')$ generates the linear system $|2M'|$.
\ep 
( This is well known if $H$ consists of all the 2-torsion points of
$B'$, for any principally polarised pair $(B',M')$.)
 
\pf 
Since the map $B'\sta{\phi_{2M'}}{\lrar}|2M'|$ is given by $a\mapsto
  t_a^*\theta +t_{-a}^*\theta$, where $\theta$ is the unique divisor in $|M'|$,
the assertion is equivalent to showing the surjectivity of the 
multiplication map
$$\sum_{\chi\in \hat{H}} H^0(M'\otimes L_\chi)\otimes H^0(M'\otimes
L_\chi)\sta{\rho}{\lrar} H^0(M'^2)..(!).$$ Here $\hat{H}$ is the dual image
of $H$ in $Pic^0(B')$. But we showed above this isomorphism, if $\hat{H}$
gives rise to a $g$-sheeted cover $(A',L')$ of $ (B',M')$, where $L'$ is of type $(1,..,2,..,4)$.
Otherwise, $\hat{H}$ gives a cover $(A',L')$ where $L'$ is of type
$(2,2,...,2)$. By similar argument used in proving (*), (!) is still true when $A'=E_1\times...\times E_g$
and $L'=L_1\times L_2...\times L_g$, where $L_j$ are line bundles of
degree 2 on the elliptic curves $E_j$. Hence our assertion is true for
a generic pair $(B',M')$. $\Box$

So, for a generic pair $(A,L)$, the map $\p H^0(L)\lrar \p H^0(M^2)$, given as
$(...,s_\chi,...)\mapsto (...,s^2_\chi,...)$ is a morphism and we obtain a
commutative diagram (I),

$
\begin{array}{ccccc}
A & \sta{\phi_L}{\lrar} & \phi_L(A) & \subset & \p^{2^g-1}=\p H^0(L)\\
\downarrow{\!\pi} & &\downarrow & &\downarrow{\!p}  \\
B=A/G & \sta{\phi_{M^2}} {\lrar} & \cK(B) & \subset & \p^{2^g-1}=\p H^0(M^2)\\
\end{array}
$

where $p(...,s_\chi,...)=(...,s_\chi^2,...)$.

\brm
Since $\phi_{M^2}\circ \pi$ is a morphism, $\phi_L$ is a morphism i.e. $L$ 
is base point free.
\erm

\bl
 Consider a pair $(A,L)$ as in 1.2. Let $\gamma \in \cG(L)\s (i)$ be an 
element of order 2. Then $H^0(L)\neq H^0(L)^{\pm}_{\gamma}$.
\el
\pf
Case 1: Suppose $\gamma=g\in \cG(L)$. Then the action of $\gamma$ is fixed
point free on $A$. Hence by Atiyah- Bott fixed point theorem,
$$h^0(L)^{+}_{\gamma}=h^0(L)^{-}_{\gamma}= h^0(L)/2.$$

Case 2: Suppose $\gamma=i$. Then 
$$h^0(L)^{\pm}_{i}= h^0(L)/2 \pm 2^{g-s-1}$$
( see [1], 4.6.6), where $s$ is the number of odd integers occurring in the
type of $L$.

Case 3: Suppose $\gamma=i.g$ and $H^0(L)=H^0(L)^+_{\gamma}$, where $g\in
\cG(L)$ is an element of order 2. Let $s\in H^0(L)^-_g$. Then $\gamma(s)=s$
gives $i(s)=-s$, i.e. $s\in H^0(L)^-_i$. Hence $H^0(L)^-_g\subset H^0(L)^-_i$
. But this contradicts the fact that $h^0(L)^-_g=2^{g-1}$ and $h^0(L)^-_i=
2^{g-1}-2^{g-s-1}$ (here $s>1$). Similarly $H^0(L)\neq H^0(L)^{-}_\gamma$.
$\Box$

Suppose $\phi_L$ is not birational and is a finite morphism of degree $d$,
$d>1$. Notice that $A\sta{\phi_{M^2}\circ \pi}{\lrar} \cK(B)$ is a Galois covering
with Galois group $(G,i)\simeq (\frac{\ze}{2\ze})^{g+1}$ and we have the extension of 
fields, $Q(\cK(B))\lrar Q(\phi_L(A))\lrar Q(A)$. Hence the Galois group of $Q(A)$
over $Q(\phi_L(A))$ is a subgroup of $(G,i)$, say $H$, of order $d$. Let
$\gamma \in H$. Then $\gamma$ is an involution on $A$, given as $a\mapsto
\epsilon a+g$ where $\epsilon = \pm 1$, $g\in G$ and it induces an involution
$\gamma^{'}$ on $H^0(L)$.

Hence $\phi_L$ factorizes as $ A\sta{\psi_1}{\lrar} A/(\gamma)\sta{\psi_2}
{\lrar}\phi_L(A)\subset \p^{2^g-1}$.
This means that the morphism $\psi_2$ is given by the pair $(N,H^0(L)^+_{
\gamma^{'}})$ or $(N^{'},H^0(L)^{-}_{\gamma^{'}})$, where $N$ and $N^{'}$ are
line bundles on $A/(\gamma)$ whose pullback to $A$ is $L$.
By 2.3, $H^0(L)\neq H^0(L)^{\pm}_{\gamma^{'}}$ and hence $\phi_L(A)$ is a 
degenerate variety in $\p^{2^g-1}$. This contradicts the fact that the morphism
$\phi_L$ is given by a complete linear system.
Hence $\phi_L$ is a birational morphism.

\section{Configuration when $g=3$}
Assume $g=3$.
Choose a $\it{theta\, structure}$ $f:\cG(L)\lrar Heis(2,4)$, ( i.e. $f$ is an
 isomorphism which restricts to identity on $\comx^*$.)
This induces an isomorphism $H^0(L)\simeq V(2,4)$ and a level
structure $K(L)\simeq \frac{\ze}{2\ze}\oplus \frac{\ze}{4\ze} \oplus \frac{\ze}{2\ze} \oplus \frac{\ze}{4\ze}$. Let $\sigma_1,
\tau_1,\sigma_2,\tau_2$ be the generators of the summands such that $o(\sigma_i)=2$ and
$o(\tau_i)=4$. The Weil form $e^L$ is given as 
$$e^L(\sigma_1,\sigma_2)=-1$$
$$e^L(\tau_1,\tau_2)=-i$$
$$e^L(\sigma_i,\tau_j)=1.$$
Then we see that the subgroup $G=<\sigma_1,\tau_1^2,\tau_2^2>$ of $K(L)$ is maximal
isotropic for the form $e^L$.

We may assume $L$ is strongly symmetric (see [10], Remark 2.4., p.160), i.e.,
 $e^L_*(g)=1 $ for all $g\in K(L)_2$,
after choosing a normalized isomorphism $\psi:L\simeq i^*(L)$,
 i.e. $\psi(0)=+1$. Here $e^L_{*}:A_2 \lrar \{\pm 1\}$ is a quadratic form
 whose value at an element $a$, of order 2 is the action of $\psi$ at the
 fibre of $L$ at $a$. 

Consider the exact sequence
$$ 1\lrar \comx^*\lrar \cG(L)\lrar K(L)\lrar 0 $$
and the homomorphism $\delta_{-1}: \cG(L) \lrar \cG(L)$, $z\mapsto
izi $.
Then $\delta_{-1}(z)=\al z^{-1}$ for some $\al\in \comx^*$.

 By [6], Proposition 2.3, p.141, we further assume that $f$ is a
 $\it{symmetric\, theta\, structure}$, i.e. $f\circ
 \delta_{-1}=D_{-1}\circ f$, where $D_{-1}:Heis(\delta)\lrar
 Heis(\delta)$ is the homomorphism $(\al,x,l)\mapsto (\al,-x,-l)$.

\bl
If $z\in \cG(L)$ is an element of order 2 and $z\neq \pm 1$ then $\delta_{-1}
(z)=e^L_*(z)z$.
\el
\pf: See [8], Proposition 3, p.309. $\Box$

\brm
Let $\sigma^{'}_1,\sigma^{'}_2,\tau^{'}_1,\tau^{'}_2\in \cG(L)$ be lifts
of $\sigma_1, \sigma_2,\tau_1,\tau_2$ such that $o(\sigma^{'}_i)=2, o(\tau^{'}
_i)=4. $
 Since $\tau_i^2\in G$, $e^L_*(\tau^2_i)=1$, hence by 3.1, $\delta_{-1}((\tau^{'}_i)^2)= (\tau^{'}_i)^2$. Hence
$\delta_{-1}(\tau^{'}_i)= c.{\tau^{'}}^{-1}_i, c=\pm 1$.
We may assume $c=+1$, by suitably altering the lift $\tau^{'}_i$.
\erm

Let $G^{'}=<\sigma^{'}_1,(\tau^{'}_1)^2,(\tau^{'}_2)^2>\subset \cG(L).$

Then $L$ descends to a principal polarization $M$ on $B=A/G$.

As remarked in Section 2,
$$H^0(L)=\oplus_{\chi\in \hat{G}}H^0(M\otimes L_\chi)$$
and $\{s_\chi\in H^0(M\otimes L_\chi),\chi \in \hat{G}\}$ form a basis of $H^0(L)$.

Consider the commutative diagram,

$
\begin{array}{ccc}
A & \sta{\psi_L}{\lrar} & Pic^0(A)  \\
\downarrow{\!\pi} &     & \uparrow{\!\hat{\pi}} \\
B & \sta{\psi_M}{\lrar} & Pic^0(B) \\
\end{array}
$

where $\psi_L(a)=t^*_aL\otimes L^{-1}$ and $\psi_M(b)=t^*_bM\otimes M^{-1}$.
Then $\psi_M$ is an isomorphism and since $\hat{\pi}(L_\chi)=0$, 
 we have $\pi^{-1}\psi^{-1}_M(L_\chi)\in K(L) \forall \chi \in \hat{G}$.
Hence $M\otimes L_\chi \simeq t^*_b M$ where $b\in \pi (K(L))$.
The basis elements $\{s_{\chi}\}_{\chi\in \hat{G}}$ can be written as
$s_0,s_1=\sigma^{'}_2(s_0),s_2=\tau^{'}_1(s_0),s_3=\tau^{'}_2(s_0),
s_4=\sigma^{'}_2\tau^{'}_1(s_0), s_5=\sigma^{'}_2\tau^{'}_2(s_0),s_6=
\tau^{'}_1 \tau^{'}_2(s_0), s_7=\sigma^{'}_2\tau^{'}_1\tau^{'}_2(s_0)$.

\bl
If $a\in K(L)_2$, then $a.i=i.a$.
\el
\pf
By 3.1, $\delta_{-1}(a)= e^L_*(a)a$. Since $e^L_*(a)=1,
a.i=i.a$. $\Box$

In particular, $g^{'}i(s_0)=ig^{'}(s_0)$, for all $ g^{'}\in G^{'}$.
Since $g^{'}s_0=s_0$, $i(s_0)\in H^0(M)$. This implies that $i(s_0)=\pm s_0$. We may assume
$i(s_0)=s_0$.

\bl
a) $i\sigma^{'}_2(s_0)= \sigma^{'}_2(s_0)$.

b) $i \tau{'}_j(s_0)=\tau^{'}_j(s_0)$.

c) $i\sigma^{'}_2\tau^{'}_j(s_0)= \sigma^{'}_2\tau^{'}_j(s_0)$.

d) $i\tau^{'}_1\tau^{'}_2(s_0)= -\tau^{'}_1\tau^{'}_2(s_0)$.

e)$i\sigma^{'}_2\tau^{'}_1\tau^{'}_2(s_0)= 
 -\sigma^{'}_2\tau^{'}_1\tau^{'}_2(s_0) $
\el
\pf
 We will use 3.3 and the fact that $g^{'}(s_0)=s_0,$ for all $g^{'}\in G^{'}$.

a) $i\sigma'_2(s_0)=\sigma'_2i(s_0)=\sigma^{'}_2(s_0)$.

b)  $i
\tau{'}_j(s_0)=\tau'^{-1}_ji(s_0)=\tau'^3_j(s_0)=\tau^{'}_j(s_0)$, (
since $\tau'^2_j\in G'$).

c)  $i\sigma^{'}_2\tau^{'}_j(s_0)= \sigma'_2i\tau'_j(s_0)=\sigma^{'}_2\tau^{'}_j(s_0)$.

d)  $i\tau^{'}_1\tau^{'}_2(s_0)= \tau'^{-1}_1i\tau'_2(s_0)=\tau'_1\tau'^2_1\tau'_2(s_0)=-\tau'_1\tau'_2\tau'^2_1(s_0)=
    -\tau^{'}_1\tau^{'}_2(s_0)$ ( since $e^L(\tau'^2_1,\tau'_2)=-1,
    \tau'^2_1\in G'$).

e) $i\sigma^{'}_2\tau^{'}_1\tau^{'}_2(s_0)= \sigma'_2i\tau^{'}_1\tau^{'}_2(s_0) 
 -\sigma^{'}_2\tau^{'}_1\tau^{'}_2(s_0) $
$\Box$

Hence we have shown the following.
\bp
The vector subspace $H^0(L)^+_i$ of $H^0(L)$ is generated by the sections
$s_0,s_1,s_2,s_3,$
$s_4,s_5 $ and the subspace $H^0(L)^{-}_i $ of $H^0(L)$ is
generated by the sections $s_6$ and $s_7$.
\ep

We then have the commutative diagram,

$
\begin{array}{ccccc}
A&\sta{\phi_L}{\lrar} & \phi_L(A) & \subset & \p (H^0(L))\\
\downarrow{\!\pi}&  &\downarrow & & \downarrow{\!p}  \\
B=A/G & \sta{\phi_{M^2}}{\lrar} & \cK(B) & \subset & \p (H^0(M^2))\\
\end{array}
$ ...(I).

Here $degree(p)=2^7$ and $degree(\pi)=8$. Since we have shown that $\phi_L$ 
is a birational morphism, $degree(\phi_L)=1$ and hence $degree(p|_{\phi_L(A)})
=2^4$. The ramification locus of $p|_{\phi_L(A)}$ is $\bigcup^7_{i=0} (H_i\cap {\phi_L(A)})$,
where $H_i$ is the hyperplane $\{s_i=0\}$ in $\p (H^0(L))$, $0\leq i\leq 7$.

Consider the group $J$ generated by the projective transformations $\al_i$,
$$(s_0,...,s_i,...,s_7)\mapsto (s_0,...,-s_i,...,s_7)$$ for $i=1,...,7$.

Then $order(J)=2^7$ and the group $J$ is the Galois group of the finite morphism
$p$.

\bp 
The group $G'\times <i>$ can be identified as a subgroup of $J$.
\ep
\pf: Since the action of $g\in G$ on the abelian threefold is fixed point free,
the $\pm1$-eigenspaces of $H^0(L)$ under the transformation $g\in G^{'}$
are equidimensional. Also, $g(s_\chi)=\chi(g).s_\chi$, for all
$\chi \in \hat{G}$, implies that $g=\al_i\al_j\al_k\al_l
\in J$, for some $0\leq i<j<k<l\leq 7$. Here $\al_0=\al_1\al_2...\al_7$.
By 3.5, $i(s_0:...:s_7)=(s_0:...:s_5:-s_6:-s_7)$. Hence the involution
$i=\al_6.\al_7$.
Hence we can identify $G^{'}\times <i>$ as a subgroup of $J$. $\Box$

Moreover, since the Galois group of the morphism $p$, $Gal(p)=J$ and
the subgroup $G'\times <i>\subset J$, leaves the image $\phi_L(A)$
invariant in $\p H^0(L)$, we have the following

\bp
Consider the commutative diagram (I). The inverse image of the variety, $\cK(B)$, has eight distinct components $h(\phi_L(A))$, where $h\in J/(G'\times <i>)$.
\ep

In Section 2, we have seen that 
$\{ t_0=s_0^2,t_1=\sigma'_2(s_0^2), t_2=\tau'_1(s_0^2),t_3=\tau'_2(s_0^2),t_4= 
\sigma'_2\tau'_1(s_0^2),t_5= \sigma'_2\tau'_2(s_0^2),t_6=\tau'_1\tau'_2(s_0^2),
t_7=\sigma'_2\tau'_1\tau'_2(s_0^2)\} $

form  a basis of $H^0(M^2)$.

\brm
( We use the same notations for the elements in $K(L)$ and their images
in $K(M^2)$.)
The elements $\sigma'_2,\tau'_1,\tau'_2$ of $\cG(M^2)$ 
 act on these sections as follows.

\begin{tabular}{llrr}
\vspace{.05in}

    & $\sigma^{'}_2$ &$ \tau^{'}_1$ & $\tau^{'}_2$  \\ 
$t_0$ &$ t_1$        &$ t_2$      &$ t_3$          \\
$t_1$ &$ t_0$        &$ t_4$      & $t_5 $       \\ 
$t_2$ &$ t_4$        &$ t_0$      &$ -t_6$         \\
$t_3$ &$ t_5$        &$ t_6$      &$t_0$          \\
$t_4$ &$ t_2$        &$ t_1$      &$-t_7$           \\
$t_5$ &$ t_3$        &$ t_7$      &$t_1$           \\
$t_6$ &$ t_7$        &$ t_3$      &$-t_2$           \\
$t_7$ &$ t_6$        &$ t_5$      &$-t_4$         \\

\end{tabular}

\erm

Now let $H_i=\{s_i=0\}$ denote the coordinate hyperplanes in $\p
H^0(L)$, for $i=0,1,...,7$.
Consider the curve $C= H_6\cap H_7\cap \phi_L(A)$. Then the involution
$i$ acts trivially on the curve $C$ and hence the degree of the
restricted morphism $\phi_L^{-1}(C)\lrar C$ is at least $2$.

\bp
The restricted morphism $\phi'_L: \phi_L^{-1}(C)\lrar C$ is of degree
$2$.
\ep
\pf: Consider the commutative diagram

$
\begin{array}{ccc}
\phi_L^{-1}(C)& \sta{\phi'_L}{\lrar} & C  \\
\downarrow{\!\pi'} &     & \downarrow{\!p'} \\
\phi_{M^2}^{-1}(p(C))& \sta{\phi'_{M^2}}{\lrar} & p(C) \\
\end{array}
$

Suppose the degree of the restricted morphism $\phi'_L$ is greater than $2$.
Since the Galois group of the morphism $\phi'_{M^2}\circ \pi'$ is the
group $G\times <i>$, the Galois group of $\phi'_L$ contains an element
$g\in G$. Hence the element $g$ acts trivially on the curve $C$. This
means that $C$ is contained in one of the eigenspaces $\p W^\pm$ of $\p
H^0(L)$, for the action of $g$. We claim that the intersection
$\phi_L(A)\cap \p W^\pm$ is at most a finite set of
points. This will give a contradiction. 

If $g^\perp=\{a\in K(L): e^L(a,g)=1\}$, then $\frac{g^\perp}{<g>}\simeq
Heis(1,1,4)$ or $ Heis(1,2,2)$ and the group $\frac{g^\perp}{<g>}$ acts
on the linear space $\p W^\pm$. Hence projecting from $\p W^\pm$ gives
a map
$\phi_g:\frac{A}{<g>}\lrar \p W^\mp$, which is base point free in the
first case ( by [2]) and has a finite base locus in the second case
( by [10]). This proves our claim.  
 $\Box$

Now, the group $G$ leaves the curve $C$ invariant and moreover
since $\sigma_2(H_6)=H_7$, we get $\sigma_2(C)=C$.
Hence the curves $$\tau_1(C)=H_3 \cap H_5\cap \phi_L(A)$$
$$ \tau_2(C)= H_2\cap H_4 \cap \phi_L(A)$$
$$ \tau_1.\tau_2(C)=H_0\cap H_1 \cap \phi_L(A)$$
are also invariant for the action of $\sigma_2$ and since for $x\in
C$, $i(x)=x$, $i.\tau_j^2(\tau_j(x))=\tau_j^2.\tau_j^{-1}i(x)=\tau_j(x)$. 
By $K(L)$-invariance of the image $\phi_L(A)$, we get

\bco
The morphism $\phi_L$ restricts to a morphism of degree $2$ on the curves $
\phi_L^{-1}(C),\,\phi_L^{-1}(\tau_1(C)),\,\phi_L^{-1}(\tau_2(C))$ and
$\phi_L^{-1}(\tau_1.\tau_2(C))$, onto their respective
images. Moreover, the Galois groups of these restricted morphisms are
$<i>,\,<i.\tau_1^2>,\,<i.\tau_2^2>$ and $<i.\tau_1^2.\tau_2^2>$,
respectively.
 
\eco

Let $A_2^+$ denote the set of points of order 2 on $A$ where the involution $i$
 acts on the fibre of $L$ at those points 
as $+1$ and $A_2^-$ denote the set of points where $i$ acts as $-1$.
By [1], Remark 4.7.7, $cardinality(A_2^+)= 48 $ and
 $cardinality(A_2^-)=16$.
 Hence if $a\in A_2^-$ and $s\in H^0(L)^+
_i$, then $s(a)=0$. This implies that for $a\in A^{-}_2$,$\phi_L(a)=
 (0:0:...:0:c_1:c_2)\in \p H^0(L)$, for some $c_1,c_2\in \comx $.

\bp
Let $a\in A^+_2$( respectively $A^{-}_2$) and $g\in K(L)_2$. 
Then $a+g\in A^+_2$( respectively $A^-_2$).
\ep
\pf: Let $g\in K(L)_2$ and $(g,\phi)\in \cG(L)$ be a lift of order $2$ and
$\psi:L\lrar i^*(L)$ be the normalized isomorphism.
By [7], Proposition 3, p.309,
$$\delta_{-1}(g,\phi)=(g,(t_g^*\psi)^{-1}\circ i^*\phi\circ\psi)$$
$$ =e^L_*(g).(g,\phi)$$
$$=(g,\phi) \, ( since \, L\, is\,strongly\,symmetric).$$
Hence the following diagram commutes

$
\begin{array}{ccc}
L& \sta{\psi}{\simeq} & i^*(L)  \\
\downarrow{\!\phi} &     & \downarrow{\!i^*(\phi)} \\
t^*_gL& \sta{t^*_g(\psi)}{\simeq} & i^*t^*_gL=t^*_g(i^*L) \\
\end{array}
$

Evaluating at $a\in  A^+_2$( respectively $A^-_2$), gives
$\psi(a)=t^*_g(\psi)(a)=\psi(a+g)$, i.e. $a+g\in  A^+_2$( respectively $A^-_2$).
$\Box$

Now let $a\in A^-_2$ then $\phi_L(a)=(0:...:c_1,c_2)$ for some
$c_1,c_2\in \comx$. Then $\sigma_2\phi_L(a)=(0:...:c_2:c_1)$. We may
assume $c_2\neq 0$.
Let $P_0=\phi_L(a)=(0:...:c:1)$ and $Q_0=p(P_0)=(0:...:c^2:1)$, for
some $c\in \comx$.

\bp
The points $h(P_0),\,h\in K(L)/<\tau^2_1,\tau^2_2>$ are of degree
$4$ on the image $\phi_L(A)$.
\ep
\pf: By 3.11, the action of $G$ on the set $A^-_2$ has two distinct
orbits, namely $O_1=\{a+g:g\in G\}$ and $O_2=\{a+\sigma_2+g: g\in
G\}$.
Then $\phi_{M^2}\circ\pi(O_1)=Q_0$ and $\phi_{M^2}\circ\pi(O_2)=
\sigma_2(Q_0)$.
Notice that $P_0\in \tau_1(C)\cap\tau_2(C)\cap
\tau_1.\tau_2(C)$. Hence, by 3.10,
$\phi^{-1}_L(P_0)=\{a,a+2\tau_1,a+2\tau_2, a+2\tau_1+2\tau_2\}$.
The assertion now follows from the $K(L)$-invariance of the image
$\phi_L(A)$.
$\Box$

\bco
The points $b(Q_0)$, where $b\in
<\pi(\sigma_2),\pi(\tau_1),\pi(\tau_2)>$, lie on the Kummer $\cK(B)$.
\eco

\section{Prym Varieties}

We recall few facts on Prym varieties ( see [5], [9], [12], for details).

Let $C$ be a smooth projective curve of genus $g$. We will assume $C$
has no vanishing theta nulls. In particular, when $g=3$, this means
$C$ is a non-hyperelliptic curve.
A point of order $2$, in $X=Jac(C)$, say $x$, defines an unramified
$2$- sheeted cover $C_x$ of $C$, $q_x:C_x\lrar C$. Let $P_x=Ker
(Nm(q_x):Jac(C_x)\lrar X)^o$, where $`o`$ denotes the connected
component containing $0\in Jac(C_x)$. Here $Nm(q_x)(\cO(\sum
r_iP_i))=\cO(\sum r_iq_x(P_i))$ is the norm map.
This defines a principally polarized abelian
variety $(P_x,\theta_{P_x})$, of dimension $g-1$.
Since the 
kernel of the dual map $q_x^{'}:X\lrar Jac(C_x)$ is generated by the
element $x$, $q_x{'}  $ induces an isomorphism $x^{\perp}/x \lrar
P_x[2]$.
Since $q_{x*}\cO_{C_x}\simeq \cO_C\oplus x$, we have $det
q_{x*}\cO_{C_x}\simeq x$.
Hence $det(q_{x*}(p)) $ is also $x$, for any $p\in ker(Nm(q_x))$.

Fix a $z\in X$ with $z^2\simeq x$. This gives a map
$$ \psi_x:Ker(Nm(q_x))\simeq P_x\cup P_x \lrar SU_C(2).$$
where $\psi_x(p)=(q_{x*}p)\otimes z$.

The image of $\psi_x$ is independent of the choice of $z$.
Recall the map
$$ SU_C(2) \sta {\phi}{\lrar} |2\theta_C|\simeq \p (H^0(SU_C(2),\cL))$$
where $\cL$ generates $Pic(SU_C(2))\simeq \ze$.

Let $\p V^+_x$ and $\p V^-_x$ be the two eigenspaces for the action of
$x$ on $|2\theta_C|$. Then there is one component of $Ker(Nm(q_x))$ in
each eigenspace. So we get a map $\phi_x:P_x\lrar \p V_x$.

\bp
The map $\phi_x:P_x\lrar \p V_x$ is the natural map
$$ P_x\lrar \cK(P_x)\subset \p (H^0(P_x,2\theta_{P_x})\simeq \p V_x.$$
\ep
\pf: See [5], Proposition 1, p.745.

\bp
For any curve $C$ and any $x$ in $X[2]- \{0\}$,
we have $ \cK(C)\cap \p V_x= \cK(P_x[2])$, ( the Schottky Jung relations).
\ep
\pf: See [5], Proposition 2 (1), p.746.

\section{Situation in $\p(H^0(L))$, when $g=3$.}

We now assume $B=J(C)$, where $J(C)$ is the Jacobian of a
non-hyperelliptic curve $C$ of genus 3.
(This is the generic situation, since the dimension of the moduli space of
principally polarized abelian threefolds is $6$ which equals the dimension of the moduli
space of curves of genus 3.)
Recall the results of Narasimhan and Ramanan ( $Theorem 1.3, Theorem
1.4$), to obtain a morphism
 $$J(C) \sta{\phi_{2\theta}}{\lrar} \cK(C)\subset F\subset
 |2\theta|$$
where 

1) $F$ is a quartic hypersurface and is the isomorphic image of
the moduli space $SU_C(2)$ and

2) the Kummer variety $\cK(C)$ is precisely the singular locus of $F$.

We will use the following
\bp
Let $L$ be an ample line bundle of type $\delta=(d_1,d_2,...,d_g)$ on
an abelian variety $A$. Then the set of irreducible representations of the
theta group $\cG(L)$, where $\al\in \comx^*$ acts as multiplication by
$\al^n$( called as of 'weight $n$'), is in bijection with the set of characters on the subgroup of
$n-$torsion elements, $K(L)_n$, of $K(L)$. Moreover, the dimension of
any such representation is $\frac{d_1.d_2...d_g}{(n,d_1)...(n,d_g)}$.
( $(n,d_i)$ denotes the greatest common divisor of $n$ and $d_i$.)
\ep
\pf: When $n=2$, the statement is proved in [6], Proposition 3.2,
p.142.
The same proof holds when $n>2$, by choosing a section over the
subgroup of $n$-torsion elements, $K(L)_n$, of $K(L)$ in the exact
sequence
$$1\lrar \comx^*\lrar \cG(L)\lrar K(L)\lrar 0$$ in the proof of [6],
Proposition 3.2.
$\Box$
\bco
The quartic $F$ in $|2\theta|$ is $\cG(2\theta)$-invariant and the
linear span of the eight cubics $\{\frac{dF}{dt_i}\}$ for
$i=0,1,...,7$ form an irreducible $\cG(2\theta)$-module where $\al\in
\comx^*$ acts as multiplication by $\al^3$.
\eco
\pf: Consider the multiplication maps $Sym^nH^0(2\theta)\sta{\rho_n}{\lrar}
H^0(2n
\theta)$. Then $I_n= Ker(\rho_n)=$ vector space of degree $n$ forms
containing the image $\cK(B)$ in $\p H^0(2\theta)$.
Since the vector spaces $Sym^nH^0(2\theta)$ and $H^0(2n\theta)$ ( via the
homomorphism $\cG(2\theta)\sta{\epsilon_n}{\lrar}\cG(2n\theta)$)
 are $\cG(2\theta)$-modules, of weight $n$ and $\rho_n$ is equivariant for
 the $\cG(2\theta)$-action, $I_n$ is also a $\cG(2\theta)$-module of weight $n$. 
Now the homogenous polynomial $F\in I_4$ and the partial derivatives
$\frac{dF}{dt_i}\in I_3$. By 5.1, it follows that $F$ is
$\cG(2\theta)$-invariant , upto scalars. If $z\in \cG(2\theta)$, then
$z\frac{dF}{dt_i}= \frac{d(zF)}{d(zt_i)}=\al\frac{dF}{d(zt_i)}\in W=
\comx\{\frac{dF}{dt_i}\}_{i=0}^7$, for some scalar $\al$. Hence $W$ is a
$\cG(2\theta)$-module of weight $3$. By 5.1, dimension of such an
irreducible representation is $8$.
This proves our assertion.
$\Box$ 

Similarly, we see that $R=F(s^2_0,...,s^2_7)$ is a $\cG(L)$-invariant
octic hypersurface in $\p H^0(L)$, by applying 5.1.
 
Recall the Weil form $e^L$ on $K(L)$ and the isotropic subgroup
$G=<\sigma_1,\tau_1^2,\tau_2^2>\subset K(L)$. Then
$e^L(\sigma_2+g,\sigma_1)=-1$, for all $g\in G$. Let $a=\sigma_2+g$,
for $g\in G$ and $a'=\sigma'_2+g'\in \cG(L)$.

Recall the basis $\{s_0,s_1,...,s_7\}$ of $H^0(L)$ and
$\{s^2_0,...,s^2_7\}$ of $H^0(M^2)$, (see Section 3).
Let $W^+_a$ and $W^-_a$ denote the eigen spaces in $H^0(L)$, for the
action of $a'$. Now $\p W_a^\pm= \{s=0: s\in W_a^\mp\}$ and $\p
V_a^+=\{t=0: t\in H^0(M^2)^-_a\}$.
Now $W^\pm _{\sigma_2}=\comx\{s_0\pm s_1,s_2\pm s_4,s_3\pm s_5,s_6\pm s_7\}$
and
$H^0(M^2)^-_{\sigma_2}=\comx\{s_0^2-s_1^2,s_2^3-s_4^2,s_3^2-s_5^2,s_6^2-s_7^2\}$.

Then $p$ restricts on $\p W^\pm _{\sigma_2}\lrar \p V^+_{\sigma_2}$ as $
(s_0;s_2:s_3,s_6)\mapsto (s_0^2:s_2^2:s_3^2:s_6^2)$, of degree $2^3$.
Similarly, one checks that if $a=\sigma_2+g,g\in G$ then $p$ restricts
to $\p W^\pm_a\lrar \p V^+_{\sigma_2}$ as 
$(z_0:...:z_3)\mapsto (z^2_0:...:z^2_3)$ of degree $2^3$.

\bp
Consider a principally polarized abelian surface $(Y,P)$, which is not
a product of elliptic curves. Let
$y_1,y_2\in Y$ be elements of order $2$, such that $e^{P^2}(y_1,y_2)=-1$. Then
we have the following.

1) There is a polarized abelian surface $(Z,N)$, such that $N$ is
   strongly symmetric of type $(1,4)$ and there is a covering map $f:Z\lrar Y$ with
   the Galois group of the map $f$ being isomorphic to $\ze/2\ze
   \times \ze/2\ze$.

2) The vector space $H^0(N)$ can be written as
$$ H^0(N)= H^0(P)\oplus H^0(t^*_{y_1}P)\oplus H^0(t^*_{y_2}P)\oplus
H^0(t^*_{y_1+y_2}P).$$

and there is a commutative diagram

$
\begin{array}{ccccc}
Z& \sta{\phi_N}{\lrar} & \phi_N(Z) & \subset & \p^3=\p H^0(N)\\
\downarrow{\!f} & &\downarrow & &\downarrow{\!q}  \\
Y& \sta{\phi_{P^2}} {\lrar} & \cK(Y) & \subset & \p^3=\p H^0(M^2)\\
\end{array}
$

 where $q(r_0:r_1:r_2:r_3)=
(r^2_0:r^2_1:r^2_2:r^2_3)$. Here $\{r_0,r_1,r_2,r_3\}$ is a basis
obtained from above decomposition of $H^0(N)$, such that $r_0,r_1,r_3\in H^0(N
)^+_i$ and $r_3\in  H^0(N)^-_i$.

\ep
\pf: 1) Consider the isomorphism $\phi_P: Y\lrar Pic^0(Y)$, $b\mapsto 
t^*_bP\otimes P^{-1}$. Let $L_{y_1}$ and $L_{y_2}$ denote the images
of $y_1$ and
$y_2$ under this map. These two line bundles define an unramified cover, 
$f:Z\lrar Y$, whose Galois group is isomorphic to $\ze/2\ze \times
\ze/2\ze$, as asserted.

 Then $N=f^*P$ is an ample line bundle and $dim H^0(N)=
4$. So to see that $N$ is of type $(1,4)$, it is enough to show that
$K(N)$ has an element of order $4$.
Consider the commutative diagram

$
\begin{array}{ccc}
Z& \sta{\psi_N}{\lrar} & Pic^0(Z)  \\
\downarrow{\!f} &     & \uparrow{\!\hat{f}} \\
Y & \sta{\psi_M}{\lrar} & Pic^0(Y) \\
\end{array}
$

Then $\hat{f}\circ \psi_M(y_i)=0$. This implies that if $z_1$ and
$z_2$ are in $Z$ such that $f(z_i)=y_i$, then $z_1,z_2\in
K(N)$. Moreover, since $e^{P^2}(y_1,y_2)=-1$ and $N^2 \simeq
f^*(P^2)$,
we have $e^{N^2}(z_1,z_2)=-1$. This gives $e^N(z_1,z_2)=\pm i$.
Hence the elements $z_1,z_2\in K(N)$ are of order $4$.

2) Clearly, $f_*N= P\oplus (P\otimes L_{y_1})\oplus (P\otimes L_{y_2})\oplus (P\otimes L_{y_1+y_2})$.
Now, in the algebraic equivalence class of $N$, there are strongly symmetric
line bundles. Hence, by tensoring $P$ with a suitable line bundle of
order $2$, we may assume that $N=f^*P$ is strongly symmetric and $r_0\in H^0(P)$ is
such that $i(r_0)=r_0$.

Since $N$ is strongly symmetric, by 3.1, $\delta_{-1}(z_j^{'})^2=(z_j^{'})^2$,
for some lifts $z_j^{'}\in \cG (N)$ of $z_j\in K(N)$.
We may further choose the lifts such that
$\delta_{-1}(z_j^{'})=(z_j^{'})^{-1}$, ( as in 3.2). In particular,
the descent data of $N$ to $P$ is $K^{'}=<(z_1^{'})^2, (z_2^{'})^2>
\subset \cG (N)$, which is a splitting over $K=<z_i^2,z_2^2>\subset K(N)$
in the exact sequence 
$$ 1\lrar \comx^* \lrar \cG (N) \lrar K(N) \lrar 0.$$
This means $(z_j^{'})^2 r_0=r_0$. Also this gives 

As in 3.5, we see that
$$ i.z_j^{'}(r_0)= z_j^{'}(r_0)$$
and $$ i.z_1^{'}.z_2^{'}(r_0)=-z_1^{'}.z_2^{'}(r_0).$$
Thus $r_0,r_1=z_1^{'}(r_0),r_2=z_2^{'}(r_0)\in H^0(N)^+_i$ and 
$r_3=z_1^{'}.z_2^{'}(r_0)\in H^0(N)^-_i$.

Hence one sees as earlier that $Gal(q)= <z^2_1,z^2_2,i>$, with a
commuatative diagram as in 5.3.
$\Box$

\bp
Let $a=\sigma_2+g,\,g\in G$ and $\p W_a$
denote an eigenspace of $a$ in $\p H^0(L)$. Then there is an abelian
surface $Z$ and a symmetric line bundle $N$ on $Z$ of type $(1,4)$ such
that $ Z\sta {\phi_{N}}{\lrar} \p (H^0(N))\sta {Heis(4)}{ \simeq} \p W_a \subset
\p H^0(L)$. Moreover, under this isomorphism, the image $\phi_N(Z)$
is mapped onto the $Heis(4)$-invariant surface $S=R\cap \p W_a$, where
$R$ is the $Heis(2,4)$- invariant hypersurface of degree 8 in $\p
H^0(L)$, defined by $F(s^2_0:s^2_1:,,,:,s^2_7)=0$. ( $F$ being the
Coble quartic).
\ep
\pf: 
Consider the restricted morphism $p:\p W_a\lrar \p V_a$, given as
$(z_0:...:z_3)\mapsto (z^2_0:...:z^2_3)$.
Then $a$ acts trivially on $\p W_a$ and $a^\perp/a(\simeq Heis(4)$) acts on $\p
W_a$, (here $a^\perp=\{y\in K(L):e^L(a,y)=1\}$). Hence there is a $Heis(4)$- action on $\p W_a$ and similarly a $Heis(2,2)$-
action on $\p V_a$. By 4.1, there is a principally polarized abelian 
surface $(P_a,\theta_{C_a})$, ($P_a$ being the Prym variety associated
to the element $\pi(a)\in K(M^2))$, such that 
$$ P_a\lrar \cK(P_a)\subset |2\theta_{C_a}|\simeq \p V_a.$$
Consider the images of $\tau_1,\tau_2$, which are elements of order 2
in $J(C)$. Since $e^{L^2}(\tau_i,a)=1$, for the
Weil form $e^{2\theta}$ on $J(C)[2]$, $\pi(\tau_1),\pi(\tau_2)\in 
\pi(a)^\perp/\pi(a)$. Moreover, $e^{2\theta}(\pi(\tau_1),\pi(\tau_2))=-1$.
  By 4.2, the points
$\phi_{M^2}\circ \pi(\tau_i)$, are nodes in the Kummer of the Prym variety
$P_a$. These nodes correspond to elements of order $2$
in $P_a$, say $\beta_1$ and $\beta_2$. Since the Weil form
$e^{2\theta_{C_a}}$ on $P_a[2]$
is induced from the Weil form $e^{2\theta}$, we have
$e^{2\theta_{C_a}}(\beta_1,\beta_2)=-1$. 
By 5.3, there is a
  polarized abelian surface $(Z,N)$ of type $(1,4)$, such that the following diagram commutes

$
\begin{array}{ccccc}
Z& \sta{\phi_N}{\lrar} & \phi_N(Z) & \subset & \p H^0(N)\\
\downarrow{\!f} & &\downarrow & &\downarrow{\!q}  \\
P_a& \sta{\phi_{2\theta_{C_a}}} {\lrar} & \cK(P_a) & \subset &|2\theta_{C_a}|\\
\end{array}
$

and for the choice of basis $\{r_0,r_1,r_2,r_3\}$, in 5.3  2),
the morphism $q$ is defined as $(r_0:r_1:r_2:r_3)\mapsto
(r^2_0:r^2_1:r^2_2:r^2_3)$, with $Gal(q)=<z^2_1,z^2_2,i>$, ($z_j$ as
in 5.3).

 Now, $R$ is the $Heis(2,4)$-invariant octic
$F(s^2_0:...:s^2_7)=0$, where $F$ is the Coble quartic.
 Note that
$S=R\cap \p W_a $ is $a^\perp/a$-invariant and is mapped onto the Kummer, $ K(P_a)$, under the restriction
morphism. Moreover, the Galois group of $p_{|S}$ is
$<\tau^2_1,\tau^2_2,i>$ which is isomorphic to the Galois group of $q$. Hence there is a $Heis(4)$- isomorphism $\p H^0(N) \lrar \p W_a$,
such that the Heisenberg invariant octic surface $\phi_N(Z)$ is mapped
onto the Heis(4)-invariant octic surface $S=R\cap \p W_a$. This proves the assertion.
$\Box$ 

It is known that the Kummer $\cK (P_a)$, has $6$ of its nodes in
each of the coordinate hyperplane, namely
the coordinate points and $3$ other distinct points.
 The preimages of the coordinate points are the coordinate points in
 $\p H^0(N)$ and $q$ is etale over the other $3$ points which are the
 pinch points of $\phi_N(Z)$ in the respective coordinate hyperplane.
\bp 
$\phi_N(Z)$ has exactly $48$ pinch points, $12$ in each coordinate
hyperplane.
\ep
\pf: See [3], Proposition 2.2, p.633.

Let $T_a$ denote the set of pinch points and the coordinate points in $\phi_N(Z
)$.

\bp
The components $h(\phi_L(A)), h\in H $ (here $H= J/(G^{'}\times i)$) and $\p W_a$ intersect at the subset $T_a$ of $ \phi_{N}(Z)$. 
In particular $\cap_{h\in H}h(\phi_L(A))= \cup _{a=\sigma_2+g, g\in
  G} T_a$. 
\ep
\pf: Since $\pi^{-1}\cK(C)= \cup_{h\in H}h(\phi_L(A))$, by 4.2 and 5.5, we
conclude that 
 $h(\phi_L(A))\cap \p
W_a=T_a$, for all $h\in H$. This gives the assertion.
$\Box$

\section{Some remarks}
a)
Consider the moduli space $\cA^l_{(1,2,4)}$ of triples $(A,c_1(L),f)$,
where $f:K(L)\lrar \ze/D\ze\times \ze/D\ze$ is a level
structure, ( here $D=(1,2,4)$).
Consider the subset of $\cA^l_{(1,2,4)}$, $\cA^{lo}_{(1,2,4)}$, parametrizing triples which
admit a $(\ze/2\ze)^3-$isogeny to the Jacobian of a non-hyperelliptic
curve.

Since $dim \cA^{lo}_{(1,2,4)}= dim \cA^l_{(1,2,4)}=6$ and $c_1(L)$
  gives a birational morphism , $\cA^{lo}_{(1,2,4)}$ is an open subset
  of $\cA^l_{(1,2,4)}$. 

Consider a triple $(A,c_1(L),f)\in \cA^{lo}_{(1,2,4)}$. We have seen
that there is a $Heis (2,4)$-invariant octic hypersurface $R$, defined
by $F(s^2_0:s^2_1:...:s^2_7)=0$, ( $F$ being the Coble quartic), such
that $\phi_L(A)\subset R\subset \p V(2,4)$.
In fact $ h(\phi_L(A))\subset Sing (R),$ for all $h\in H$, ( $H$ as in 5.6).

Now $F$ is a $Heis(2,2,2)$-invariant quartic polynomial in $\p V(2,2,2)$.
Since the space of $Heis(2,2,2)$-invariant quartics is
$14$-dimensional, ( see [4], p.186]), the space of
$Heis(2,4)$-invariant octics in $\p^7$ which are of the form
$R=F(s_0^2:...:s_7^2)$ where $F$ is a $Heis(2,2,2)$-invariant quartic, is
also $14$-dimensional. Call this space as
 
$P(Sym^8V(2,4)^{Heis(2,4)'})=\p^{14}$.

So there is a morphism   
$$ \cA^{lo}_{(1,2,4)}\sta {T}{\lrar}  \p^{14}$$
where $T$ is defined as $ (A,c_1(L),f)\mapsto R$.

One may try to study this morphism, from a moduli point of view.

b) Consider the special basis $\{s_0^2,...,s_7^2\}$ ( which is
different from the usual $\it{Heisenberg}$ basis) of $H^0(2\theta)$ and the action 
of the elements of the subgroup $<\sigma_2,\tau_1^2,\tau_2^2>\subset K(2\theta)$ on this
basis ( see 3.8).

Also, by 3.12, the points
$b(P_0)\in \phi_L(A)$,
where $b\in <\sigma_2,\tau_1,\tau_2>\subset K(L)$, $P_0=(0:...:0:c:1)$
and the point $Q_0=(0:...:0:c^2:1)\in \cK(C)$, for some non-zero
$c\in \comx$. 
With these data, in addition to knowing the geometry of $SU_C(2)$ in
$|2\theta|$- linear system  one may try to know the equation of the $\it{Coble
  \,quartic}$, in terms of this basis $\{s_0^2,...,s_7^2\}$.

 \begin{thebibliography}{99}
\bib [1]{1}  Birkenhake, Ch., Lange, H. : {\em Complex abelian varieties},   
 Springer-Verlag, Berlin, (1992).

\bib [2]{2} Birkenhake, Ch.; Lange, H.; Ramanan, S.{ Primitive line
    bundles on abelian threefolds}, Manuscripta Math.$\bf{81}$,
  no. 3-4, 299-310 (1993).

\bib [3]{3} Birkenhake, Ch., Lange, H., van Straten, D. :{\em Abelian surfaces
of type (1,4)}, Math.Annalen.$\bf{285}$, 625-646 (1989).

\bib [4]{4} Dolgachev, I., Ortland, D. :{ \em Point sets in
  projective spaces }, Asterisque No. 165, (1988).

\bib [5]{5} B.van Geemen, Previato,E. :{\em Prym varieties and the
  Verlinde formula}, Math.Annalen.$\bf{294}$, 741-754, (1992).

\bib [6]{6} Iyer, J. :{\em Projective normality of abelian surfaces, given by
primitive line bundles }, Manuscr.Math.$\bf{98}$, 139-153, (1999).

\bib [7]{7} Kouvidakis, A., Pantev,T. :{\em The automorphism group of
  the moduli space of semistable vector bundles },
  Math.Annalen.$\bf{302}$, no 2, 225-268 (1995).

\bib [8]{8}  Mumford, D. : {\em On the equations defining Abelian
varieties 1}, Invent. math. $\bf{1}$, 287-354 (1966). 

\bib [9]{9}  Mumford, D. :{\em Prym varieties. 1.} In : Contributions
  to analysis, London, New York: Academic Press, 325-350, 1974.

\bib [10]{10}  Nagaraj, D,S., Ramanan, S.: {\em Polarisations of type
$(1,2,...,2)$ on abelian varieties}, Duke Math.Journal, $\bf{ 80}$, No.1,
157-194, (1995). 

\bib [11]{11}  Narasimhan, M.S., Ramanan, S. : {\em Moduli of vector
  bundles on a compact Riemann surface}, Ann. of Math. (2) $\bf{89}$,
  14-51, (1969).

\bib [12]{12} Narasimhan, M.S., Ramanan, S. : {\em $2\theta$-Linear Systems},
 Vector bundles on algebraic varieties (Bombay, 1984) 415-427, Tata Inst.Fund.
Res.Stud.Math.11, Tata Inst.Fund.Res, Bombay, 1987

\end {thebibliography}

\end{document}